\documentclass[12pt]{article}
\usepackage{eurosym}
\usepackage{amsfonts}
\usepackage{hyperref}
\usepackage{color}
\usepackage{amsmath,amssymb,latexsym}
\usepackage{graphicx}

\newtheorem{theorem}{Theorem}[section]
\newtheorem{lemma}[theorem]{Lemma}
\newtheorem{definition}[theorem]{Definition}

\newtheorem{example}[theorem]{Example}

\numberwithin{equation}{section}

\begin{document}

\title{{\Huge \textbf{On a generalization of fractional Langevin equation}}}
\author{Saeed Kosari$\sp{a}$, Milad Yadollahzadeh$\sp{b}$%
\thanks{Corresponding author. m.yadollahzadeh@yahoo.com}, Zehui Shao$\sp{a}$
~and Yongsheng Rao$\sp{a}$ \\
$\sp{a}${\small \textit{Institute of Computing Science and Technology, }%
}\\
{\small \textit{\ Guangzhou university, Guangzhou, China}}\\
$\sp{b}${\small \textit{Department of Mathematics, Faculty of Mathematical
Sciences, }}\\
{\small \textit{University of Mazandaran, Babolsar, Iran}}\\
}
\date{}
\maketitle

\begin{abstract}
In this work, we consider a generalization of the nonlinear Langevin
equation of fractional orders with boundary value conditions. The existence
and uniqueness of solutions are studied by using results of the fixed point
theory. Moreover, the previous results of fractional Langevin equations are
a special case of our problem.\newline


\end{abstract}

\section{Introduction}

The Langevin equation has been applied to describe phenomena that have
stochastic properties. First, Langevin equation \cite{Langevin} used for a
particle submerged in a fluid with Brownian motion. Langevin equation has
the various versions which proposed in different fields. For instance,
Pineda and Stamatakis \cite{Pineda} investigated the stochastic modelling of
surface reactions by using reflected chemical Langevin equations. In \cite%
{lu}, based on Langevin equation, molecular dynamics simulation was used in
the study of structural, thermal properties of matter in different phases.
Moreover, the modified Langevin equation was applied as a macroscopic
stochastic nonlinear model of many geophysical processes in \cite{Czechowski}%
. Mendoza-M\'{e}ndez \textit{et al.} \cite{Mendoza} proposed a dynamic
equivalence between the longtime dynamic properties of atomic and colloidal
liquids by a more formal fundamental derivation of the generalized Langevin
equation for a tracer particle in an atomic liquid.

Study of existence results of equations with applications in engineering
industries and modeling different processes is gaining much importance and
attention. Recently, Fallah and Mehrdoust studied \cite{Fallah} the
existence and uniqueness of the solution to the stochastic differential
equation of the double Heston model which is defined by two independent
variance processes with non-Lipschitz diffusion.
The existence of global mild solutions was obtained for the incompressible
nematic liquid crystal flow in the whole space in \cite{Zhang}. Marasi \textit{et al.} \cite{Marasi}
studied the existence and uniqueness of positive solutions for nonlinear fractional differential equation boundary value problems by using new fixed point results.
For some recent works on existence
results, see \cite{Boussandel,Chen,Lv,Pavlova,Tan,Tang}

According to the applications of Langevin equation, on the existence theory
of solutions for this class of equations is an important subject. Ahmad and
Nieto \cite{Nieto} discussed the existence of solutions of nonlinear
Langevin equation involving two fractional orders with Dirichlet boundary
conditions as follows:
\begin{eqnarray*}
&&^{c}D_{0+}^{\beta }\left( ^{c}D_{0+}^{\alpha }+\lambda \right) x\left(
t\right) =f\left( t,x(t)\right) ,\text{ }t\in \left( 0,1\right) , \\
&&x\left( 0\right) =\gamma _{1},\text{ }x\left( 1\right) =\gamma _{2},\ \
0<\alpha ,\beta \leqslant 1,
\end{eqnarray*}%
where $^{c}D^{\alpha }$ is the Caputo fractional derivative of order $\alpha
$, $f:[0,1]\times \mathbb{X}\rightarrow \mathbb{X}$ is a continuous
function, $\mathbb{X}$ is a Banach space, $\lambda $ is a real number and $%
\gamma _{1},\gamma _{2}\in \mathbb{X}$. Also, they  introduce a $q$-fractional variant of nonlinear Langevin equation of different orders with $q$-fractional antiperiodic boundary conditions in \cite{Ahmad}.

Fazli and J. Nieto \cite{Fazli} investigated the existence and uniqueness of
solutions for nonlinear Langevin equation of fractional orders with
anti-periodic boundary conditions as:%
\begin{eqnarray*}
&&^{c}D_{0+}^{\beta }\left( ^{c}D_{0+}^{\alpha }+\lambda \right) x\left(
t\right) =f\left( t,x(t)\right) ,\text{ }t\in \left( 0,1\right) ,\text{ }%
0<\alpha \leqslant 1,\text{ }1<\beta \leqslant 2, \\
&&x\left( 0\right) +x\left( 1\right) =0,\text{ \ }^{c}D_{0+}^{\alpha
}x\left( 0\right) +\text{ }^{c}D_{0+}^{\alpha }x\left( 1\right) =0,\text{ \ }%
\mathcal{D}_{0+}^{2\alpha }x\left( 0\right) +\mathcal{D}_{0+}^{2\alpha
}x\left( 1\right) =0,
\end{eqnarray*}%
where the function $f:[0,1]\times \mathbb{R}\rightarrow \mathbb{R}$ is
continuous and $\lambda $ is a real number and $\mathcal{D}_{0+}^{2\alpha }$
is the sequential fractional derivative.

This manuscript is concerned to study the existence and uniqueness of
solutions a generalization of the nonlinear Langevin equation of different
fractional orders with four-point boundary conditions provided as:
\begin{eqnarray}
&&^{c}D_{0+}^{\beta }\left( ^{c}D_{0+}^{\alpha }+\gamma \right) x\left(
t\right) =f\left( t,x(t),x^{\prime }\left( t\right) \right) ,\text{ }t\in
\left( 0,1\right) ,  \label{e1} \\
&&x\left( 0\right) =x(1)=x^{\prime }\left( 0\right) =x^{\prime }\left(
1\right) =0,\ \ 0<\alpha \leqslant 1,\text{ }2<\beta \leqslant 3,  \label{e2}
\end{eqnarray}%
where $^{c}D^{\alpha }$ is the Caputo fractional derivative of order $\alpha
$, $f:[0,1]\times \mathbb{R\times R}\rightarrow \mathbb{R}$ is a continuous
function and $\gamma $ is a real number. Under suitable assumptions on part
of nonlinear function and by application of the fixed point theory, we
obtain the existence results of the considered problem.

The paper is scheduled as follows: In Section 2 are recalled some necessary
preliminaries from fractional calculus. Next, we obtain the Green functions
corresponding to the problem. In section 3 is devoted to the existence and
uniqueness of solutions for boundary value problem (\ref{e1}) and (\ref{e2}%
). Also, an example to illustrate our results is given. Finally, we propose
some conclusions.

\section{Preliminaries}

In this section, we present some definitions and lemmas which are needed for
our results.

\begin{definition}
\cite{Kilbas,Pod} The Riemann-Liouville fractional integral of order $\rho
>0 $ of a function $g:(0,\infty )\rightarrow \mathbb{R}$ is defined by
\begin{equation*}
I_{0^{+}}^{\rho }g(\omega )=\frac{1}{\Gamma (\rho )}\int_{0}^{\omega
}(\omega -s)^{\rho -1}g(s)ds,
\end{equation*}%
provided the right-hand side is pointwise defined on $(0,\infty ).$
\end{definition}

\begin{definition}
\cite{Kilbas,Pod} The Riemann-Liouville fractional derivative of order $\rho
>0$ of a continuous function $g:(0,\infty )\rightarrow \mathbb{R}$ is
defined by
\begin{equation*}
D_{0^{+}}^{\rho }g(\omega )=\frac{1}{\Gamma (m-\rho )}\left( \frac{d}{%
d\omega }\right) ^{m}\int_{0}^{\omega }(\omega -s)^{m-\rho -1}g(s)ds,
\end{equation*}%
where $m=[\rho ]+1,$ provided that right-hand side is pointwise defined on $%
(0,\infty ).$
\end{definition}

\begin{definition}
\cite{Kilbas,Pod} For a function $g$ given on the interval $[0,\infty ),$
the Caputo fractional derivative of order $\rho >0$ of $g$ is defined by
\begin{equation*}
^{c}D_{0^{+}}^{\rho }g(\omega )=\frac{1}{\Gamma (m-\rho )}\int_{0}^{\omega
}(\omega -s)^{m-\alpha -1}g^{\left( n\right) }(s)ds,
\end{equation*}%
where $m=[\rho ]+1.$
\end{definition}

\begin{lemma}
\label{proper}\cite{Kilbas,Pod} Let $\alpha ,\beta \geqslant 0$ and $n\in
\mathbb{N}$, then the following relations hold:

\begin{enumerate}
\item $D_{a+}^{\alpha }I_{a+}^{\alpha }f(t)=f(t)$,

\item $D_{a+}^{\beta }I_{a+}^{\alpha }f(t)=D_{a+}^{\beta -\alpha }f(t),$ \
(if $\beta \geqslant \alpha $),

\item $D_{a+}^{\beta }I_{a+}^{\alpha }f(t)=I_{a+}^{\alpha -\beta }f(t),$ \
(if $\alpha \geqslant \beta $),

\item $I_{a+}^{\alpha }\left( t-a\right) ^{\beta }=\frac{\Gamma \left( \beta
+1\right) }{\Gamma \left( \alpha +\beta +1\right) }\left( t-a\right)
^{\alpha +\beta },$

\item $D_{a+}^{\alpha }\left( t-a\right) ^{\beta }=\frac{\Gamma \left( \beta
+1\right) }{\Gamma \left( \beta -\alpha +1\right) }\left( t-a\right) ^{\beta
-\alpha }$.
\end{enumerate}
\end{lemma}

\begin{lemma}
\cite{Kilbas} Let $\rho >0$. Then, the fractional differential equation $%
^{c}D^{\rho }u(t)=0,$ has a general solution as
\begin{equation*}
u(t)=a_{0}+a_{1}t+a_{2}t^{2}+\cdots +a_{m-1}t^{m-1},\ \
\end{equation*}%
for some $a_{k}\in \mathbb{R}$, $k=1,\ldots ,m-1$, and $m=[\rho ]+1.$
\end{lemma}

\begin{lemma}
\cite{Kilbas}\label{lem2} Let $\rho >0$. Then, we have
\begin{equation*}
I_{0^{+}}^{\rho }~^{c}D_{0^{+}}^{\rho }u(t)=u(t)+a_{0}+a_{1}t+\cdots
+a_{m-1}t^{m-1},
\end{equation*}%
where $a_{k}\in \mathbb{R}$, $k=1,\ldots ,m-1$, and $m=[\rho ]+1.$
\end{lemma}

\begin{lemma}
\label{lem3}Let $y\in C[0,1],$ a unique solution of boundary value problem
for following fractional Langevin equation%
\begin{eqnarray*}
&&^{c}D_{0+}^{\beta }\left( ^{c}D_{0+}^{\alpha }+\gamma \right) x\left(
t\right) =y\left( t\right) ,\text{ }t\in \left( 0,1\right) , \\
&&x\left( 0\right) =x(1)=x^{\prime }\left( 0\right) =x^{\prime }\left(
1\right) =0,\ \ 0<\alpha \leqslant 1,\text{ }2<\beta \leqslant 3,
\end{eqnarray*}%
is given by%
\begin{equation*}
x(t)=\int_{0}^{1}G\left( t,s\right) y(s)ds+\int_{0}^{1}H\left( t,s\right)
x(s)ds,
\end{equation*}%
where%
\begin{equation*}
G\left( t,s\right) =\frac{1}{\Gamma \left( \alpha +\beta \right) }\left\{
\begin{array}{ll}
\begin{array}{l}
\left( t-s\right) ^{\alpha +\beta -1}+\left[ \left( \alpha +1\right)
t^{\alpha +2}-\left( \alpha +2\right) t^{\alpha +1}\right] \left( 1-s\right)
^{\alpha +\beta -1} \\
+\left( \alpha +\beta -1\right) \left( t^{\alpha +1}-t^{\alpha +2}\right)
\left( 1-s\right) ^{\alpha +\beta -2},%
\end{array}
& s\leqslant t, \\
&  \\
\begin{array}{l}
\left[ \left( \alpha +1\right) t^{\alpha +2}-\left( \alpha +2\right)
t^{\alpha +1}\right] \left( 1-s\right) ^{\alpha +\beta -1} \\
+\left( \alpha +\beta -1\right) \left( t^{\alpha +1}-t^{\alpha +2}\right)
\left( 1-s\right) ^{\alpha +\beta -2},%
\end{array}
& t\leqslant s,%
\end{array}%
\right.
\end{equation*}%
and%
\begin{equation*}
H\left( t,s\right) =\frac{1}{\Gamma \left( \alpha \right) }\left\{
\begin{array}{ll}
\begin{array}{l}
-\gamma \left( t-s\right) ^{\alpha +\beta -1}+\gamma \left[ \left( \alpha
+2\right) t^{\alpha +1}-\left( \alpha +1\right) t^{\alpha +2}\right] \left(
1-s\right) ^{\alpha -1} \\
+\gamma \left( \alpha -1\right) \left( t^{\alpha +2}-t^{\alpha +1}\right)
\left( 1-s\right) ^{\alpha -2},%
\end{array}
& s\leqslant t, \\
&  \\
\begin{array}{l}
\gamma \left[ \left( \alpha +2\right) t^{\alpha +1}-\left( \alpha +1\right)
t^{\alpha +2}\right] \left( 1-s\right) ^{\alpha -1} \\
+\gamma \left( \alpha -1\right) \left( t^{\alpha +2}-t^{\alpha +1}\right)
\left( 1-s\right) ^{\alpha -2},%
\end{array}
& t\leqslant s.%
\end{array}%
\right.
\end{equation*}
\end{lemma}

\noindent \textbf{Proof.} For $2<\beta \leqslant 3$, Lemma \ref{lem2} yields%
\begin{equation*}
^{c}D^{\alpha }x\left( t\right) =I_{0^{+}}^{\beta }y(t)-\gamma x\left(
t\right) +c_{0}+c_{1}t+c_{2}t^{2},
\end{equation*}%
where $c_{0},c_{1},c_{2}\in \mathbb{R}.$ Also, for $0<\alpha \leqslant 1,$
we have%
\begin{equation}
x\left( t\right) =I_{0^{+}}^{\alpha +\beta }y(t)-\gamma I_{0^{+}}^{\alpha
}x\left( t\right) +c_{0}\frac{t^{\alpha }}{\Gamma \left( \alpha +1\right) }%
+c_{1}\frac{t^{\alpha +1}}{\Gamma \left( \alpha +2\right) }+c_{2}\frac{%
2t^{\alpha +2}}{\Gamma \left( \alpha +3\right) }+c_{3}.  \label{a0}
\end{equation}%
By the condition $x\left( 0\right) =0,$ we give $c_{3}=0.$ Differentiation
of (\ref{a0}) with respect to $t$ produces
\begin{equation}
x^{\prime }\left( t\right) =I_{0^{+}}^{\alpha +\beta -1}y(t)-\gamma
I_{0^{+}}^{\alpha -1}x\left( t\right) +c_{0}\frac{t^{\alpha -1}}{\Gamma
\left( \alpha \right) }+c_{1}\frac{t^{\alpha }}{\Gamma \left( \alpha
+1\right) }+c_{2}\frac{2t^{\alpha +1}}{\Gamma \left( \alpha +2\right) }.
\label{z1}
\end{equation}%
The condition $x^{\prime }\left( 0\right) =0$, for (\ref{z1}), implies that $%
c_{0}=0.$ Now, by conditions $x\left( 1\right) =0$ and $x^{\prime }\left(
1\right) =0,$ we give%
\begin{equation}
\frac{1}{\Gamma \left( \alpha +\beta \right) }\int_{0}^{1}\left( 1-s\right)
^{\alpha +\beta -1}y(s)ds-\frac{\gamma }{\Gamma \left( \alpha \right) }%
\int_{0}^{1}\left( 1-s\right) ^{\alpha -1}x(s)ds+c_{1}\frac{1}{\Gamma \left(
\alpha +2\right) }+c_{2}\frac{2}{\Gamma \left( \alpha +3\right) }=0,
\label{q1}
\end{equation}%
and%
\begin{equation}
\frac{1}{\Gamma \left( \alpha +\beta -1\right) }\int_{0}^{1}\left(
1-s\right) ^{\alpha +\beta -2}y(s)ds-\frac{\gamma }{\Gamma \left( \alpha
-1\right) }\int_{0}^{1}\left( 1-s\right) ^{\alpha -2}x(s)ds+c_{1}\frac{1}{%
\Gamma \left( \alpha +1\right) }+c_{2}\frac{2}{\Gamma \left( \alpha
+2\right) }=0,  \label{q2}
\end{equation}%
respectively. By solving the system of (\ref{q1}) and (\ref{q2}), we obtain
\begin{eqnarray*}
c_{1} &=&\frac{\Gamma \left( \alpha +2\right) }{\Gamma \left( \alpha +\beta
-1\right) }\int_{0}^{1}\left( 1-s\right) ^{\alpha +\beta -2}y(s)ds-\frac{%
\Gamma \left( \alpha +3\right) }{\Gamma \left( \alpha +\beta \right) }%
\int_{0}^{1}\left( 1-s\right) ^{\alpha +\beta -1}y(s)ds \\
&&+\frac{\gamma \Gamma \left( \alpha +3\right) }{\Gamma \left( \alpha
\right) }\int_{0}^{1}\left( 1-s\right) ^{\alpha -1}x(s)ds-\frac{\gamma
\Gamma \left( \alpha +2\right) }{\Gamma \left( \alpha -1\right) }%
\int_{0}^{1}\left( 1-s\right) ^{\alpha -2}x(s)ds,
\end{eqnarray*}%
and%
\begin{eqnarray*}
c_{2} &=&\frac{\left( \alpha +1\right) \Gamma \left( \alpha +3\right) }{%
2\Gamma \left( \alpha +\beta \right) }\int_{0}^{1}\left( 1-s\right) ^{\alpha
+\beta -1}y(s)ds-\frac{\Gamma \left( \alpha +3\right) }{2\Gamma \left(
\alpha +\beta -1\right) }\int_{0}^{1}\left( 1-s\right) ^{\alpha +\beta
-2}y(s)ds \\
&&+\frac{\gamma \Gamma \left( \alpha +3\right) }{2\Gamma \left( \alpha
-1\right) }\int_{0}^{1}\left( 1-s\right) ^{\alpha -2}x(s)ds-\frac{\gamma
\left( \alpha +1\right) \Gamma \left( \alpha +3\right) }{2\Gamma \left(
\alpha \right) }\int_{0}^{1}\left( 1-s\right) ^{\alpha -1}x(s)ds.
\end{eqnarray*}%
By substituting the values of $c_{0},c_{1},c_{2},c_{2}$ in (\ref{a0}), then
we give
\begin{eqnarray*}
&&x(t)\underset{}{=}\frac{1}{\Gamma \left( \alpha +\beta \right) }%
\int_{0}^{t}\left( t-s\right) ^{\alpha +\beta -1}y(s)ds-\frac{\gamma }{%
\Gamma \left( \alpha \right) }\int_{0}^{t}\left( t-s\right) ^{\alpha
-1}x(s)ds \\
&&\underset{}{+}\left[ \left( \alpha +2\right) t^{\alpha +1}-\left( \alpha
+1\right) t^{\alpha +2}\right] \frac{\gamma }{\Gamma \left( \alpha \right) }%
\int_{0}^{1}\left( 1-s\right) ^{\alpha -1}x(s)ds+\frac{\gamma \left(
t^{\alpha +2}-t^{\alpha +1}\right) }{\Gamma \left( \alpha -1\right) }%
\int_{0}^{1}\left( 1-s\right) ^{\alpha -2}x(s)ds \\
&&+\frac{t^{\alpha +1}-t^{\alpha +2}}{\Gamma \left( \alpha +\beta -1\right) }%
\int_{0}^{1}\left( 1-s\right) ^{\alpha +\beta -2}y(s)ds+\frac{\left( \alpha
+1\right) t^{\alpha +2}-\left( \alpha +2\right) t^{\alpha +1}}{\Gamma \left(
\alpha +\beta \right) }\int_{0}^{1}\left( 1-s\right) ^{\alpha +\beta
-1}y(s)ds \\
&&\underset{}{=}\int_{0}^{t}\left( \frac{\left( t-s\right) ^{\alpha +\beta
-1}+\left[ \left( \alpha +1\right) t^{\alpha +2}-\left( \alpha +2\right)
t^{\alpha +1}\right] \left( 1-s\right) ^{\alpha +\beta -1}}{\Gamma \left(
\alpha +\beta \right) }\right. \\
&&\left. +\frac{\left( \alpha +\beta -1\right) \left( t^{\alpha
+1}-t^{\alpha +2}\right) \left( 1-s\right) ^{\alpha +\beta -2}}{\Gamma
\left( \alpha +\beta \right) }\right) y(s)ds \\
&&+\int_{t}^{1}\frac{\left[ \left( \alpha +1\right) t^{\alpha +2}-\left(
\alpha +2\right) t^{\alpha +1}\right] \left( 1-s\right) ^{\alpha +\beta
-1}+\left( \alpha +\beta -1\right) \left( t^{\alpha +1}-t^{\alpha +2}\right)
\left( 1-s\right) ^{\alpha +\beta -2}}{\Gamma \left( \alpha +\beta \right) }%
y(s)ds \\
&&+\int_{0}^{t}\left( \frac{-\gamma \left( t-s\right) ^{\alpha +\beta
-1}+\gamma \left[ \left( \alpha +2\right) t^{\alpha +1}-\left( \alpha
+1\right) t^{\alpha +2}\right] \left( 1-s\right) ^{\alpha -1}}{\Gamma \left(
\alpha \right) }\right. \\
&&\left. +\frac{\gamma \left( \alpha -1\right) \left( t^{\alpha
+2}-t^{\alpha +1}\right) \left( 1-s\right) ^{\alpha -2}}{\Gamma \left(
\alpha \right) }\right) x(s)ds \\
&&+\int_{t}^{1}\frac{\gamma \left[ \left( \alpha +2\right) t^{\alpha
+1}-\left( \alpha +1\right) t^{\alpha +2}\right] \left( 1-s\right) ^{\alpha
-1}+\gamma \left( \alpha -1\right) \left( t^{\alpha +2}-t^{\alpha +1}\right)
\left( 1-s\right) ^{\alpha -2}}{\Gamma \left( \alpha \right) }x(s)ds \\
&&\underset{}{=}\int_{0}^{1}G\left( t,s\right) y(s)ds+\int_{0}^{1}H\left(
t,s\right) x(s)ds.\text{ \ }\square
\end{eqnarray*}

\section{Main results}

In this section, we propose the existence results of Langevin differential
equation of fractional orders (\ref{e1})-(\ref{e2}). Let the Banach space $%
B=C^{1}[0,1]$ be equipped with the norm:
\begin{equation}
\Vert x\Vert =\max_{t\in \lbrack 0,1]}|x(t)|+\max_{t\in \lbrack
0,1]}|x^{\prime }(t)|.  \label{Norm Banach}
\end{equation}

To prove the main results, we need the following assumptions:

\begin{itemize}
\item[$A_{1})$] $f:\left[ 0,1\right] \times \mathbb{R}\times \mathbb{R}%
\rightarrow \mathbb{R}$ is continuous function.

\item[$A_{2})$] There exists a constant $w>0$ such that
\begin{equation*}
|f(t,x,y)-f(t,u,v)|\leq w\left( |x-u|+|y-v|\right) ,
\end{equation*}%
for all $t\in \lbrack 0,1],$ $x,y,u,v\in \mathbb{R}.$

\item[$A_{3})$] There exists a nonnegative function $\sigma \in L[0,1]$ such
that\newline
\begin{equation*}
|f(t,x,y)|\leq \sigma (t)+a_{1}|x|^{\tau _{1}}+a_{2}|y|^{\tau _{2}},
\end{equation*}%
where $a_{1},a_{2}\in \mathbb{R}_{+}$ and $0<\tau _{1},\tau _{2}<1$.
\end{itemize}

For the sake of convenience, we define the following constants:

$\mathcal{K}_{1}=\left\{ \frac{4\left( \alpha +1\right) +2\beta }{\Gamma
\left( \alpha +\beta +1\right) }\right\} ,\mathcal{K}_{2}=\left\{ \frac{%
2\alpha +4}{\Gamma \left( \alpha +\beta \right) }+\frac{2\left( \alpha
+1\right) \left( \alpha +2\right) }{\Gamma \left( \alpha +\beta +1\right) }%
\right\} ,\mathcal{K}=\mathcal{K}_{1}+\mathcal{K}_{2},$

$\mathcal{L}_{1}=\left\{ \frac{\left\vert \gamma \right\vert 4\left( \alpha
+1\right) }{\Gamma \left( \alpha +1\right) }\right\} ,\mathcal{L}%
_{2}=\left\{ \frac{\left\vert \gamma \right\vert \left( 2\alpha +4\right) }{%
\Gamma \left( \alpha \right) }+\frac{2\left\vert \gamma \right\vert \left(
\alpha +1\right) \left( \alpha +2\right) }{\Gamma \left( \alpha +1\right) }%
\right\},\mathcal{L}=\mathcal{L}_{1}+\mathcal{L}_{2}.$

\begin{theorem}
\label{Th1}Under the hypotheses $\left( A_{1}\right) $ and $\left(
A_{3}\right) $, the fractional Langevin equation (\ref{e1})-(\ref{e2}) has a
solution.
\end{theorem}

\noindent \textbf{Proof.} We define the operator $\mathcal{Z}:B\rightarrow B$
as follow:%
\begin{equation*}
\left( \mathcal{Z}x\right) (t)=\int_{0}^{1}G\left( t,s\right) f\left(
s,x(s),x^{\prime }\left( s\right) \right) ds+\int_{0}^{1}H\left( t,s\right)
x(s)ds.
\end{equation*}%
Lemma \ref{lem3} implies that the fixed points of the operator $\mathcal{Z}$
are the same solutions of the boundary value problem (\ref{e1})-(\ref{e2}).
We consider a ball $U_{r}=\left\{ x\in B,\Vert x\Vert \leq r\right\} $ so
that\newline
$\max \left\{ 4\mathcal{K}\Vert \sigma \Vert ,\left( 4a_{1}\mathcal{K}%
\right) ^{1/1-\tau _{1}},\left( 4a_{2}\mathcal{K}\right) ^{1/1-\tau _{2}},4r%
\mathcal{L}\right\} \leq r.$ For any $x\in U_{r}$ and by $\left(
A_{3}\right) $, we show that $\mathcal{Z}U_{r}\subset U_{r}$, then%
\begin{eqnarray*}
&&\left\vert \left( \mathcal{Z}x\right) (t)\right\vert \underset{}{\leq }%
\frac{1}{\Gamma \left( \alpha +\beta \right) }\int_{0}^{t}\left( t-s\right)
^{\alpha +\beta -1}\left\vert f\left( s,x(s),x^{\prime }\left( s\right)
\right) \right\vert ds \\
&&+\frac{t^{\alpha +1}+t^{\alpha +2}}{\Gamma \left( \alpha +\beta -1\right) }%
\int_{0}^{1}\left( 1-s\right) ^{\alpha +\beta -2}\left\vert f\left(
s,x(s),x^{\prime }\left( s\right) \right) \right\vert ds \\
&&+\frac{\left( \alpha +1\right) t^{\alpha +2}+\left( \alpha +2\right)
t^{\alpha +1}}{\Gamma \left( \alpha +\beta \right) }\int_{0}^{1}\left(
1-s\right) ^{\alpha +\beta -1}\left\vert f\left( s,x(s),x^{\prime }\left(
s\right) \right) \right\vert ds+\frac{\left\vert \gamma \right\vert }{\Gamma
\left( \alpha \right) }\int_{0}^{t}\left( t-s\right) ^{\alpha -1}\left\vert
x(s)\right\vert ds
\end{eqnarray*}
\begin{eqnarray*}
&&+\left[ \left( \alpha +2\right) t^{\alpha +1}+\left( \alpha +1\right)
t^{\alpha +2}\right] \frac{\left\vert \gamma \right\vert }{\Gamma \left(
\alpha \right) }\int_{0}^{1}\left( 1-s\right) ^{\alpha -1}\left\vert
x(s)\right\vert ds+\frac{\left\vert \gamma \right\vert \left( t^{\alpha
+2}+t^{\alpha +1}\right) }{\Gamma \left( \alpha -1\right) }%
\int_{0}^{1}\left( 1-s\right) ^{\alpha -2}\left\vert x(s)\right\vert ds \\
&&\underset{}{\leq }\frac{1}{\Gamma \left( \alpha +\beta \right) }%
\int_{0}^{t}\left( t-s\right) ^{\alpha +\beta -1}\left( \sigma
(s)+a_{1}|x(s)|^{\tau _{1}}+a_{2}|x^{\prime }\left( s\right) |^{\tau
_{2}}\right) ds \\
&&+\frac{2}{\Gamma \left( \alpha +\beta -1\right) }\int_{0}^{1}\left(
1-s\right) ^{\alpha +\beta -2}\left( \sigma (s)+a_{1}|x(s)|^{\tau
_{1}}+a_{2}|x^{\prime }\left( s\right) |^{\tau _{2}}\right) ds \\
&&+\frac{2\alpha +3}{\Gamma \left( \alpha +\beta \right) }\int_{0}^{1}\left(
1-s\right) ^{\alpha +\beta -1}\left( \sigma (s)+a_{1}|x(s)|^{\tau
_{1}}+a_{2}|x^{\prime }\left( s\right) |^{\tau _{2}}\right) ds \\
&&+\frac{\left\vert \gamma \right\vert \left\Vert x\right\Vert }{\Gamma
\left( \alpha \right) }\int_{0}^{t}\left( t-s\right) ^{\alpha -1}ds+\frac{%
\left( 2\alpha +3\right) \left\vert \gamma \right\vert \left\Vert
x\right\Vert }{\Gamma \left( \alpha \right) }\int_{0}^{1}\left( 1-s\right)
^{\alpha -1}ds+\frac{2\left\vert \gamma \right\vert \left\Vert x\right\Vert
}{\Gamma \left( \alpha -1\right) }\int_{0}^{1}\left( 1-s\right) ^{\alpha
-2}ds \\
&&\underset{}{\leq }\left( \left\Vert \sigma \right\Vert +a_{1}r^{\tau
_{1}}+a_{2}r^{\tau _{2}}\right) \left\{ \frac{t^{\alpha +\beta }}{\Gamma
\left( \alpha +\beta +1\right) }+\frac{2}{\Gamma \left( \alpha +\beta
\right) }+\frac{2\alpha +3}{\Gamma \left( \alpha +\beta +1\right) }\right\}
\\
&&+\left\Vert x\right\Vert \left\{ \frac{\left\vert \gamma \right\vert
t^{\alpha }}{\Gamma \left( \alpha +1\right) }+\frac{\left( 2\alpha +3\right)
\left\vert \gamma \right\vert }{\Gamma \left( \alpha +1\right) }+\frac{%
2\left\vert \gamma \right\vert }{\Gamma \left( \alpha \right) }\right\} \\
&&\underset{}{\leq }\left( \left\Vert \sigma \right\Vert +a_{1}r^{\tau
_{1}}+a_{2}r^{\tau _{2}}\right) \left\{ \frac{4\left( \alpha +1\right)
+2\beta }{\Gamma \left( \alpha +\beta +1\right) }\right\} +r\left\{ \frac{%
\left\vert \gamma \right\vert 4\left( \alpha +1\right) }{\Gamma \left(
\alpha +1\right) }\right\} \\
&&\underset{}{=}\left( \left\Vert \sigma \right\Vert +a_{1}r^{\tau
_{1}}+a_{2}r^{\tau _{2}}\right) \mathcal{K}_{1}+r\mathcal{L}_{1}
\end{eqnarray*}%
and%
\begin{eqnarray*}
&&\left\vert \left( \mathcal{Z}x\right) ^{\prime }\left( t\right)
\right\vert \underset{}{\leq }\frac{1}{\Gamma \left( \alpha +\beta -1\right)
}\int_{0}^{t}\left( t-s\right) ^{\alpha +\beta -2}\left\vert f\left(
s,x(s),x^{\prime }\left( s\right) \right) \right\vert ds \\
&&+\frac{\left( \alpha +1\right) t^{\alpha }+\left( \alpha +2\right)
t^{\alpha +1}}{\Gamma \left( \alpha +\beta -1\right) }\int_{0}^{1}\left(
1-s\right) ^{\alpha +\beta -2}\left\vert f\left( s,x(s),x^{\prime }\left(
s\right) \right) \right\vert ds \\
&&+\frac{\left( \alpha +1\right) \left( \alpha +2\right) \left( t^{\alpha
+1}+t^{\alpha }\right) }{\Gamma \left( \alpha +\beta \right) }%
\int_{0}^{1}\left( 1-s\right) ^{\alpha +\beta -1}\left\vert f\left(
s,x(s),x^{\prime }\left( s\right) \right) \right\vert ds \\
&&+\frac{\left\vert \gamma \right\vert }{\Gamma \left( \alpha -1\right) }%
\int_{0}^{t}\left( t-s\right) ^{\alpha -2}\left\vert x(s)\right\vert ds+%
\frac{\left\vert \gamma \right\vert \left( \alpha +1\right) \left( \alpha
+2\right) \left( t^{\alpha }+t^{\alpha +1}\right) }{\Gamma \left( \alpha
\right) }\int_{0}^{1}\left( 1-s\right) ^{\alpha -1}\left\vert
x(s)\right\vert ds \\
&&+\frac{\left\vert \gamma \right\vert \left[ \left( \alpha +2\right)
t^{\alpha +1}+\left( \alpha +1\right) t^{\alpha }\right] }{\Gamma \left(
\alpha -1\right) }\int_{0}^{1}\left( 1-s\right) ^{\alpha -2}\left\vert
x(s)\right\vert ds \\
&&\underset{}{\leq }\left( \left\Vert \sigma \right\Vert +a_{1}r^{\tau
_{1}}+a_{2}r^{\tau _{2}}\right) \left\{ \frac{t^{\alpha +\beta }}{\Gamma
\left( \alpha +\beta \right) }+\frac{2\alpha +3}{\Gamma \left( \alpha +\beta
\right) }+\frac{2\left( \alpha +1\right) \left( \alpha +2\right) }{\Gamma
\left( \alpha +\beta +1\right) }\right\} \\
&&+\left\Vert x\right\Vert \left\{ \frac{\left\vert \gamma \right\vert
t^{\alpha -1}}{\Gamma \left( \alpha \right) }+\frac{2\left\vert \gamma
\right\vert \left( \alpha +1\right) \left( \alpha +2\right) }{\Gamma \left(
\alpha +1\right) }+\frac{\left\vert \gamma \right\vert \left( 2\alpha
+3\right) }{\Gamma \left( \alpha \right) }\right\} \\
&&\underset{}{\leq }\left( \left\Vert \sigma \right\Vert +a_{1}r^{\tau
_{1}}+a_{2}r^{\tau _{2}}\right) \left\{ \frac{2\alpha +4}{\Gamma \left(
\alpha +\beta \right) }+\frac{2\left( \alpha +1\right) \left( \alpha
+2\right) }{\Gamma \left( \alpha +\beta +1\right) }\right\} \\
&&+r\left\{ \frac{\left\vert \gamma \right\vert \left( 2\alpha +4\right) }{%
\Gamma \left( \alpha \right) }+\frac{2\left\vert \gamma \right\vert \left(
\alpha +1\right) \left( \alpha +2\right) }{\Gamma \left( \alpha +1\right) }%
\right\} \\
&&\underset{}{=}\left( \left\Vert \sigma \right\Vert +a_{1}r^{\tau
_{1}}+a_{2}r^{\tau _{2}}\right) \mathcal{K}_{2}+r\mathcal{L}_{2}.
\end{eqnarray*}%
So%
\begin{eqnarray*}
\left\Vert \mathcal{Z}x\right\Vert &=&\max \left\vert \left( \mathcal{Z}%
x\right) \left( t\right) \right\vert +\max \left\vert \left( \mathcal{Z}%
x\right) ^{\prime }\left( t\right) \right\vert \\
&\leq &\left( \left\Vert \sigma \right\Vert +a_{1}r^{\tau _{1}}+a_{2}r^{\tau
_{2}}\right) \left( \mathcal{K}_{1}+\mathcal{K}_{2}\right) +r\left( \mathcal{%
L}_{1}+\mathcal{L}_{2}\right) \\
&=&\left( \left\Vert \sigma \right\Vert +a_{1}r^{\tau _{1}}+a_{2}r^{\tau
_{2}}\right) \mathcal{K}+r\mathcal{L} \\
&\leq &\frac{r}{4}+\frac{r}{4}+\frac{r}{4}+\frac{r}{4}=r.
\end{eqnarray*}%
Next, we prove that the operator $\mathcal{Z}$ is completely continuous. The
functions $f$, $G\left( t,s\right) $ and $H\left( t,s\right) $ are
continuous, hence the operator $\mathcal{Z}$ is continuous. Let $M=\underset{%
t\in \lbrack 0,1],x\in U_{r}}{\max }|f\left( t,x(t),x^{\prime }\left(
t\right) \right) |+1,$ for any $x\in U_{r}$ and $t_{1},t_{2}\in \lbrack 0,1]$
such that $t_{1}<t_{2},$ we have%
\begin{eqnarray}
&&\left\vert \left( \mathcal{Z}x\right) \left( t_{2}\right) -\left( \mathcal{%
Z}x\right) \left( t_{1}\right) \right\vert \underset{}{\leq }\left\vert
\frac{1}{\Gamma \left( \alpha +\beta \right) }\int_{0}^{t_{2}}\left(
t_{2}-s\right) ^{\alpha +\beta -1}f\left( s,x(s),x^{\prime }\left( s\right)
\right) ds\right.  \notag \\
&&\left. -\frac{1}{\Gamma \left( \alpha +\beta \right) }\int_{0}^{t_{1}}%
\left( t_{1}-s\right) ^{\alpha +\beta -1}f\left( s,x(s),x^{\prime }\left(
s\right) \right) ds\right\vert  \notag \\
&&+\left\vert \frac{\left( t_{2}^{\alpha +1}-t_{1}^{\alpha +1}\right)
+\left( t_{2}^{\alpha +2}-t_{1}^{\alpha +2}\right) }{\Gamma \left( \alpha
+\beta -1\right) }\int_{0}^{1}\left( 1-s\right) ^{\alpha +\beta -2}f\left(
s,x(s),x^{\prime }\left( s\right) \right) ds\right\vert  \notag \\
&&+\left\vert \frac{\left( \alpha +1\right) \left( t_{2}^{\alpha
+2}-t_{1}^{\alpha +2}\right) +\left( \alpha +2\right) \left( t_{2}^{\alpha
+1}-t_{1}^{\alpha +1}\right) }{\Gamma \left( \alpha +\beta \right) }%
\int_{0}^{1}\left( 1-s\right) ^{\alpha +\beta -1}f\left( s,x(s),x^{\prime
}\left( s\right) \right) ds\right\vert  \notag \\
&&+\left\vert \frac{\gamma }{\Gamma \left( \alpha \right) }%
\int_{0}^{t_{2}}\left( t_{2}-s\right) ^{\alpha -1}x(s)ds-\frac{\gamma }{%
\Gamma \left( \alpha \right) }\int_{0}^{t_{1}}\left( t_{1}-s\right) ^{\alpha
-1}x(s)ds\right\vert  \notag \\
&&+\left\vert \left[ \left( \alpha +2\right) \left( t_{2}^{\alpha
+1}-t_{1}^{\alpha +1}\right) +\left( \alpha +1\right) \left( t_{2}^{\alpha
+2}-t_{1}^{\alpha +2}\right) \right] \frac{\gamma }{\Gamma \left( \alpha
\right) }\int_{0}^{1}\left( 1-s\right) ^{\alpha -1}x(s)ds\right\vert  \notag
\\
&&+\left\vert \gamma \frac{\left( t_{2}^{\alpha +1}-t_{1}^{\alpha +1}\right)
+\left( t_{2}^{\alpha +2}-t_{1}^{\alpha +2}\right) }{\Gamma \left( \alpha
-1\right) }\int_{0}^{1}\left( 1-s\right) ^{\alpha -2}x(s)ds\right\vert
\notag \\
&&\underset{}{\leq }\frac{2M}{\Gamma \left( \alpha +\beta +1\right) }\left(
t_{2}-t_{1}\right) ^{\alpha +\beta }+\frac{M}{\Gamma \left( \alpha +\beta
+1\right) }\left( t_{2}^{\alpha +\beta }-t_{1}^{\alpha +\beta }\right)
\notag \\
&&+M\frac{\left( t_{2}^{\alpha +1}-t_{1}^{\alpha +1}\right) +\left(
t_{2}^{\alpha +2}-t_{1}^{\alpha +2}\right) }{\Gamma \left( \alpha +\beta
\right) }+M\frac{\left( \alpha +1\right) \left( t_{2}^{\alpha
+2}-t_{1}^{\alpha +2}\right) +\left( \alpha +2\right) \left( t_{2}^{\alpha
+1}-t_{1}^{\alpha +1}\right) }{\Gamma \left( \alpha +\beta +1\right) }
\notag \\
&&+\frac{2r\left\vert \gamma \right\vert }{\Gamma \left( \alpha +1\right) }%
\left( t_{2}-t_{1}\right) ^{\alpha }+\frac{r\left\vert \gamma \right\vert }{%
\Gamma \left( \alpha +1\right) }\left( t_{2}^{\alpha }-t_{1}^{\alpha
}\right) +\frac{r\left\vert \gamma \right\vert }{\Gamma \left( \alpha
\right) }\left[ \left( t_{2}^{\alpha +1}-t_{1}^{\alpha +1}\right) +\left(
t_{2}^{\alpha +2}-t_{1}^{\alpha +2}\right) \right]  \notag \\
&&+\frac{r\left\vert \gamma \right\vert }{\Gamma \left( \alpha +1\right) }%
\left[ \left( \alpha +2\right) \left( t_{2}^{\alpha +1}-t_{1}^{\alpha
+1}\right) +\left( \alpha +1\right) \left( t_{2}^{\alpha +2}-t_{1}^{\alpha
+2}\right) \right] .  \label{d11}
\end{eqnarray}%
Also,%
\begin{eqnarray}
&&\left\vert \left( \mathcal{Z}x\right) ^{\prime }\left( t_{2}\right)
-\left( \mathcal{Z}x\right) ^{\prime }\left( t_{1}\right) \right\vert
\underset{}{\leq }\left\vert \frac{1}{\Gamma \left( \alpha +\beta -1\right) }%
\int_{0}^{t_{2}}\left( t_{2}-s\right) ^{\alpha +\beta -2}f\left(
s,x(s),x^{\prime }\left( s\right) \right) ds\right.  \notag \\
&&\left. -\frac{1}{\Gamma \left( \alpha +\beta -1\right) }%
\int_{0}^{t_{1}}\left( t_{1}-s\right) ^{\alpha +\beta -2}f\left(
s,x(s),x^{\prime }\left( s\right) \right) ds\right\vert  \notag \\
&&+\left\vert \frac{\left( \alpha +1\right) \left( t_{2}^{\alpha
}-t_{1}^{\alpha }\right) +\left( \alpha +2\right) \left( t_{2}^{\alpha
+1}-t_{1}^{\alpha +1}\right) }{\Gamma \left( \alpha +\beta -1\right) }%
\int_{0}^{1}\left( 1-s\right) ^{\alpha +\beta -2}f\left( s,x(s),x^{\prime
}\left( s\right) \right) ds\right\vert  \notag \\
&&+\left\vert \frac{\left( \alpha +1\right) \left( \alpha +2\right) \left[
\left( t_{2}^{\alpha +1}-t_{1}^{\alpha +1}\right) +\left( t_{2}^{\alpha
}-t_{1}^{\alpha }\right) \right] }{\Gamma \left( \alpha +\beta \right) }%
\int_{0}^{1}\left( 1-s\right) ^{\alpha +\beta -1}f\left( s,x(s),x^{\prime
}\left( s\right) \right) ds\right\vert  \notag \\
&&+\left\vert \frac{\gamma }{\Gamma \left( \alpha -1\right) }%
\int_{0}^{t_{2}}\left( t_{2}-s\right) ^{\alpha -2}x(s)ds-\frac{\gamma }{%
\Gamma \left( \alpha -1\right) }\int_{0}^{t_{1}}\left( t_{1}-s\right)
^{\alpha -2}x(s)ds\right\vert  \notag \\
&&+\left\vert \frac{\gamma \left( \alpha +1\right) \left( \alpha +2\right) %
\left[ \left( t_{2}^{\alpha +1}-t_{1}^{\alpha +1}\right) +\left(
t_{2}^{\alpha }-t_{1}^{\alpha }\right) \right] }{\Gamma \left( \alpha
\right) }\int_{0}^{1}\left( 1-s\right) ^{\alpha -1}x(s)ds\right\vert  \notag
\\
&&+\left\vert \frac{\gamma \left[ \left( \alpha +2\right) \left(
t_{2}^{\alpha +1}-t_{1}^{\alpha +1}\right) +\left( \alpha +1\right) \left(
t_{2}^{\alpha }-t_{1}^{\alpha }\right) \right] }{\Gamma \left( \alpha
-1\right) }\int_{0}^{1}\left( 1-s\right) ^{\alpha -2}x(s)ds\right\vert
\notag \\
&&\underset{}{\leq }\frac{2M}{\Gamma \left( \alpha +\beta \right) }\left(
t_{2}-t_{1}\right) ^{\alpha +\beta -1}+\frac{M}{\Gamma \left( \alpha +\beta
\right) }\left( t_{2}^{\alpha +\beta -1}-t_{1}^{\alpha +\beta -1}\right)
\notag \\
&&+\frac{M}{\Gamma \left( \alpha +\beta \right) }\left[ \left( \alpha
+1\right) \left( t_{2}^{\alpha }-t_{1}^{\alpha }\right) +\left( \alpha
+2\right) \left( t_{2}^{\alpha +1}-t_{1}^{\alpha +1}\right) \right]  \notag
\\
&&+\frac{M}{\Gamma \left( \alpha +\beta +1\right) }\left( \alpha +1\right)
\left( \alpha +2\right) \left[ \left( t_{2}^{\alpha +1}-t_{1}^{\alpha
+1}\right) +\left( t_{2}^{\alpha }-t_{1}^{\alpha }\right) \right]  \notag \\
&&+\frac{2r\left\vert \gamma \right\vert }{\Gamma \left( \alpha \right) }%
\left( t_{2}-t_{1}\right) ^{\alpha -1}+\frac{r\left\vert \gamma \right\vert
}{\Gamma \left( \alpha \right) }\left( t_{2}^{\alpha -1}-t_{1}^{\alpha
-1}\right) +\frac{r\left\vert \gamma \right\vert \left( \alpha +1\right)
\left( \alpha +2\right) }{\Gamma \left( \alpha +1\right) }\left[ \left(
t_{2}^{\alpha +1}-t_{1}^{\alpha +1}\right) +\left( t_{2}^{\alpha
}-t_{1}^{\alpha }\right) \right]  \notag \\
&&+\frac{r\left\vert \gamma \right\vert }{\Gamma \left( \alpha \right) }%
\left[ \left( \alpha +2\right) \left( t_{2}^{\alpha +1}-t_{1}^{\alpha
+1}\right) +\left( \alpha +1\right) \left( t_{2}^{\alpha }-t_{1}^{\alpha
}\right) \right] .  \label{d2}
\end{eqnarray}%
By (\ref{d11}) and (\ref{d2}), clearly that the functions $\left(
t_{2}-t_{1}\right) ^{\alpha +\beta -i},$ $t_{2}^{\alpha +\beta
-i}-t_{1}^{\alpha +\beta -i},$ $\left( t_{2}-t_{1}\right) ^{\alpha -i}$ ($%
i=0,1$) and $t_{2}^{\alpha +j}-t_{1}^{\alpha +j}$ ($j=0,1,2$) are uniformly
continuous on $[0,1]$. Then, $\mathcal{Z}\left( U_{r}\right) $ is
equicontinuous and the Arzela--Ascoli theorem implies that $\overline{%
\mathcal{Z}\left( U_{r}\right) }$ is compact, hence the operator $\mathcal{Z}%
:U_{r}\rightarrow U_{r}$ is completely continuous. Therefore, by the
Schauder fixed-point theorem, we conclude that the problem (\ref{e1}) and (%
\ref{e2}) has a solution. $\square $

By applying Banach fixed point theorem, we prove the uniqueness of solution
of the problem (\ref{e1}) and (\ref{e2}).

\begin{theorem}
Let the assumptions $(H_{1}$-$H_{3})$ are satisfied, then the boundary value
problem (\ref{e1}) and (\ref{e2}) has a uniqueness solution provided that $%
\psi =\psi _{1}+\psi _{2}<1,$ where%
\begin{equation*}
\psi _{1}=\left\{ \frac{14w}{\Gamma \left( \alpha +\beta +1\right) }+\frac{%
8\left\vert \gamma \right\vert }{\Gamma \left( \alpha +1\right) }\right\} ,%
\text{ }\psi _{2}=\left\{ \frac{36w}{\Gamma \left( \alpha +\beta +1\right) }+%
\frac{18\left\vert \gamma \right\vert }{\Gamma \left( \alpha +1\right) }%
\right\} .
\end{equation*}
\end{theorem}

\noindent \textbf{Proof.} For any $x,y\in B,$ $t\in \left[ 0,1\right] $ and
by condition $(H_{2})$, we give%
\begin{eqnarray}
&&\left\vert \left( \mathcal{Z}x\right) \left( t\right) -\left( \mathcal{Z}%
y\right) \left( t\right) \right\vert \underset{}{\leq }\frac{1}{\Gamma
\left( \alpha +\beta \right) }\int_{0}^{t}\left( t-s\right) ^{\alpha +\beta
-1}\left\vert f\left( s,x(s),x^{\prime }(s)\right) -f\left( s,y(s),y^{\prime
}(s)\right) \right\vert ds  \notag \\
&&+\frac{\left( \alpha +1\right) t^{\alpha +2}+\left( \alpha +2\right)
t^{\alpha +1}}{\Gamma \left( \alpha +\beta \right) }\int_{0}^{1}\left(
1-s\right) ^{\alpha +\beta -1}\left\vert f\left( s,x(s),x^{\prime
}(s)\right) -f\left( s,y(s),y^{\prime }(s)\right) \right\vert ds  \notag \\
&&+\frac{t^{\alpha +1}+t^{\alpha +2}}{\Gamma \left( \alpha +\beta -1\right) }%
\int_{0}^{1}\left( 1-s\right) ^{\alpha +\beta -2}\left\vert f\left(
s,x(s),x^{\prime }(s)\right) -f\left( s,y(s),y^{\prime }(s)\right)
\right\vert ds  \notag \\
&&+\frac{\left\vert \gamma \right\vert }{\Gamma \left( \alpha \right) }%
\int_{0}^{t}\left( t-s\right) ^{\alpha -1}\left\vert x(s)-y(s)\right\vert ds
\notag \\
&&+\frac{\left\vert \gamma \right\vert \left[ \left( \alpha +2\right)
t^{\alpha +1}+\left( \alpha +1\right) t^{\alpha +2}\right] }{\Gamma \left(
\alpha \right) }\int_{0}^{1}\left( 1-s\right) ^{\alpha -1}\left\vert
x(s)-y(s)\right\vert ds  \notag \\
&&+\frac{\left\vert \gamma \right\vert \left[ t^{\alpha +2}+t^{\alpha +1}%
\right] }{\Gamma \left( \alpha -1\right) }\int_{0}^{1}\left( 1-s\right)
^{\alpha -2}\left\vert x(s)-y(s)\right\vert ds  \notag \\
&&\underset{}{\leq }\frac{w\left\Vert x-y\right\Vert }{\Gamma \left( \alpha
+\beta \right) }\int_{0}^{t}\left( t-s\right) ^{\alpha +\beta -1}ds+\frac{%
5w\left\Vert x-y\right\Vert }{\Gamma \left( \alpha +\beta \right) }%
\int_{0}^{1}\left( 1-s\right) ^{\alpha +\beta -1}ds  \notag \\
&&+\frac{2w\left\Vert x-y\right\Vert }{\Gamma \left( \alpha +\beta -1\right)
}\int_{0}^{1}\left( 1-s\right) ^{\alpha +\beta -2}ds+\frac{\left\vert \gamma
\right\vert \left\Vert x-y\right\Vert }{\Gamma \left( \alpha \right) }%
\int_{0}^{t}\left( t-s\right) ^{\alpha -1}ds  \notag \\
&&+\frac{5\left\vert \gamma \right\vert \left\Vert x-y\right\Vert }{\Gamma
\left( \alpha \right) }\int_{0}^{1}\left( 1-s\right) ^{\alpha -1}ds+\frac{%
2\left\vert \gamma \right\vert \left\Vert x-y\right\Vert }{\Gamma \left(
\alpha -1\right) }\int_{0}^{1}\left( 1-s\right) ^{\alpha -2}ds  \notag \\
&&\underset{}{\leq }w\left\Vert x-y\right\Vert \left\{ \frac{t^{\alpha
+\beta }}{\Gamma \left( \alpha +\beta +1\right) }+\frac{5}{\Gamma \left(
\alpha +\beta +1\right) }+\frac{2}{\Gamma \left( \alpha +\beta \right) }%
\right\}  \notag \\
&&+\left\Vert x-y\right\Vert \left\{ \frac{\left\vert \gamma \right\vert
t^{\alpha }}{\Gamma \left( \alpha +1\right) }+\frac{5\left\vert \gamma
\right\vert }{\Gamma \left( \alpha +1\right) }+\frac{2\left\vert \gamma
\right\vert }{\Gamma \left( \alpha \right) }\right\}  \notag \\
&&\underset{}{\leq }\left\Vert x-y\right\Vert \left\{ \frac{14w}{\Gamma
\left( \alpha +\beta +1\right) }+\frac{8\left\vert \gamma \right\vert }{%
\Gamma \left( \alpha +1\right) }\right\} \underset{}{=}\psi _{1}\left\Vert
x-y\right\Vert ,  \label{s1}
\end{eqnarray}%
and%
\begin{eqnarray}
&&\left\vert \left( \mathcal{Z}x\right) ^{\prime }\left( t\right) -\left(
\mathcal{T}y\right) ^{\prime }\left( t\right) \right\vert \underset{}{\leq }%
\frac{1}{\Gamma \left( \alpha +\beta -1\right) }\int_{0}^{t}\left(
t-s\right) ^{\alpha +\beta -2}\left\vert f\left( s,x(s),x^{\prime
}(s)\right) -f\left( s,y(s),y^{\prime }(s)\right) \right\vert ds  \notag \\
&&+\frac{\left( \alpha +1\right) \left( \alpha +2\right) \left( t^{\alpha
+1}+t^{\alpha }\right) }{\Gamma \left( \alpha +\beta \right) }%
\int_{0}^{1}\left( 1-s\right) ^{\alpha +\beta -1}\left\vert f\left(
s,x(s),x^{\prime }(s)\right) -f\left( s,y(s),y^{\prime }(s)\right)
\right\vert ds  \notag \\
&&+\frac{\left( \alpha +1\right) t^{\alpha }+\left( \alpha +2\right)
t^{\alpha +1}}{\Gamma \left( \alpha +\beta -1\right) }\int_{0}^{1}\left(
1-s\right) ^{\alpha +\beta -2}\left\vert f\left( s,x(s),x^{\prime
}(s)\right) -f\left( s,y(s),y^{\prime }(s)\right) \right\vert ds  \notag \\
&&+\frac{\left\vert \gamma \right\vert }{\Gamma \left( \alpha -1\right) }%
\int_{0}^{t}\left( t-s\right) ^{\alpha -2}\left\vert x(s)-y(s)\right\vert ds
\notag \\
&&+\frac{\left\vert \gamma \right\vert \left( \alpha +1\right) \left( \alpha
+2\right) \left( t^{\alpha }+t^{\alpha +1}\right) }{\Gamma \left( \alpha
\right) }\int_{0}^{1}\left( 1-s\right) ^{\alpha -1}\left\vert
x(s)-y(s)\right\vert ds  \notag \\
&&+\frac{\left\vert \gamma \right\vert \left[ \left( \alpha +2\right)
t^{\alpha +1}+\left( \alpha +1\right) t^{\alpha }\right] }{\Gamma \left(
\alpha -1\right) }\int_{0}^{1}\left( 1-s\right) ^{\alpha -2}\left\vert
x(s)-y(s)\right\vert ds  \notag
\end{eqnarray}
\begin{eqnarray}
&&\underset{}{\leq }\frac{w\left\Vert x-y\right\Vert }{\Gamma \left( \alpha
+\beta -1\right) }\int_{0}^{t}\left( t-s\right) ^{\alpha +\beta -2}ds+\frac{%
12w\left\Vert x-y\right\Vert }{\Gamma \left( \alpha +\beta \right) }%
\int_{0}^{1}\left( 1-s\right) ^{\alpha +\beta -1}ds  \notag \\
&&+\frac{5w\left\Vert x-y\right\Vert }{\Gamma \left( \alpha +\beta -1\right)
}\int_{0}^{1}\left( 1-s\right) ^{\alpha +\beta -2}ds+\frac{\left\vert \gamma
\right\vert \left\Vert x-y\right\Vert }{\Gamma \left( \alpha -1\right) }%
\int_{0}^{t}\left( t-s\right) ^{\alpha -2}ds  \notag \\
&&+\frac{12\left\vert \gamma \right\vert \left\Vert x-y\right\Vert }{\Gamma
\left( \alpha \right) }\int_{0}^{1}\left( 1-s\right) ^{\alpha -1}ds+\frac{%
5\left\vert \gamma \right\vert \left\Vert x-y\right\Vert }{\Gamma \left(
\alpha -1\right) }\int_{0}^{1}\left( 1-s\right) ^{\alpha -2}ds  \notag \\
&&\underset{}{\leq }w\left\Vert x-y\right\Vert \left\{ \frac{t^{\alpha
+\beta -1}}{\Gamma \left( \alpha +\beta \right) }+\frac{12}{\Gamma \left(
\alpha +\beta +1\right) }+\frac{5}{\Gamma \left( \alpha +\beta \right) }%
\right\}  \notag \\
&&+\left\Vert x-y\right\Vert \left\{ \frac{\left\vert \gamma \right\vert
t^{\alpha -1}}{\Gamma \left( \alpha \right) }+\frac{12\left\vert \gamma
\right\vert }{\Gamma \left( \alpha +1\right) }+\frac{5\left\vert \gamma
\right\vert }{\Gamma \left( \alpha \right) }\right\}  \notag \\
&&\underset{}{\leq }\left\Vert x-y\right\Vert \left\{ \frac{36w}{\Gamma
\left( \alpha +\beta +1\right) }+\frac{18\left\vert \gamma \right\vert }{%
\Gamma \left( \alpha +1\right) }\right\} \underset{}{=}\psi _{2}\left\Vert
x-y\right\Vert .  \label{s22}
\end{eqnarray}%
Using (\ref{s1}) and (\ref{s22}), we obtain
\begin{equation*}
\left\Vert \mathcal{Z}x-\mathcal{Z}y\right\Vert \leq \psi \left\Vert
x-y\right\Vert ,
\end{equation*}%
where $\psi <1$. Hence the operator $\mathcal{Z}$ is a contraction operator
and the contraction mapping principle implies that the problem (\ref{e1})
and (\ref{e2}) has a unique solution. $\square $

Here, we investigate the existence of solutions to a special of fractional
Langevin equation in following example.

\begin{example}
\label{exp1} Consider the following nonlinear Langevin equation of
fractional orders%
\begin{equation}
\left\{
\begin{array}{l}
^{c}D_{0+}^{\frac{5}{2}}\left( ^{c}D_{0+}^{\frac{1}{2}}+\frac{1}{5}\right)
u\left( t\right) =f(t,u(t),u^{\prime }\left( t\right) ),\text{ }t\in \left(
0,1\right) , \\
u\left( 0\right) =u\left( 1\right) =0,\text{ \ \ }u^{\prime }\left( 0\right)
=u^{\prime }\left( 1\right) =0,%
\end{array}%
\right.  \label{ex1}
\end{equation}%
where $f(t,u(t),u^{\prime }(t))=\tan ^{-1}t+\left( t-\frac{1}{3}\right)
^{2}(u(t))^{\tau _{1}}+\frac{t}{e}(u^{\prime }(t))^{\tau _{2}}$, $0<\tau
_{1},\tau _{2}<1$. Observe that the function $f$ is continuous, also
\begin{eqnarray*}
\left\vert f(t,u(t),u^{\prime }\left( t\right) )\right\vert &\leq &t+(t-%
\frac{1}{3})^{2}|u(t)|^{\tau _{1}}+\frac{t}{e}|u^{\prime }(t)|^{\tau _{2}} \\
&\leq &1+\frac{4}{9}|u(t)|^{\tau _{1}}+\frac{1}{e}|u^{\prime }(t)|^{\tau
_{2}}.
\end{eqnarray*}%
Thus, the assumptions (H$_{1}$) and (H$_{2}$) are satisfied and Theorem \ref%
{Th1} implies that the problem (\ref{ex1}) has a solution.
\end{example}

\section{Conclusions}

We have provided a generalization of the nonlinear Langevin equation of
fractional orders with boundary conditions. By considering the Banach space $%
C^{1}[0,1]$ be equipped with the norm $\Vert x\Vert =\max_{t\in \lbrack
0,1]}|x(t)|+\max_{t\in \lbrack 0,1]}|x^{\prime }(t)|$, under the assumptions
$\left( H_{1}\right) $ and $\left( H_{3}\right) $, we have shown the
existence of solution for BVP (\ref{e1}) and (\ref{e2}). Then, the
uniqueness of solution of the problem is proved by Banach fixed point
theorem. Also, we have proposed an example to illustrate the results.

\end{document}